\documentclass{amsart}
\usepackage{amssymb, array, amscd}

\theoremstyle{plain}
\newtheorem{theorem}{Theorem}[section]
\newtheorem{corollary}[theorem]{Corollary}
\newtheorem{lemma}[theorem]{Lemma}
\newtheorem{proposition}[theorem]{Proposition}
\newtheorem{algorithm}[theorem]{Algorithm}
\newtheorem{definition}[theorem]{Definition}
\newtheorem{example}[theorem]{Example}

\newtheorem{remark}[theorem]{Remark}
\numberwithin{equation}{section}

\DeclareMathOperator{\Ker}{Ker}
\DeclareMathOperator{\ord}{ord}
\DeclareMathOperator{\ex}{exp}
\DeclareMathOperator{\length}{length}
\newcommand{\la}{\langle}
\newcommand{\ra}{\rangle}
\newcommand{\N}{\mathbb{N}}
\newcommand{\Z}{\mathbb{Z}}
\newcommand{\R}{\mathcal{R}}

\newcommand{\G}{\mathcal{G}}
\newcommand{\x}{\mathbf{x}}
\newcommand{\bfc}{\mathbf{c}}

\begin{document}

\title[On the additive structure of indecomposable modules]{Computing the additive structure of  indecomposable modules over Dedekind-like rings using Gr\"obner bases.}
\author{Maria A. Avi\~no-Diaz and Luis D. Garcia-Puente}
\address{Department of Mathematics, UPR-Cayey , PR 00736}
\email{mavino@uprr.pr}
\address{ Department of Mathematics\ Texas A\&M University\\
College Station, TX 77843-3368}
\email{lgp@math.tamu.edu}
\thanks{The first author is supported by the National
Institute of Health,  PROGRAM SCORE, 2004-08, 546112, University of Puerto Rico-Rio Piedras Campus,
IDEA Network of Biomedical Research Excellence, and the Laboratory Gauss University of Puerto Rico Research.
She wants to thank Professor O. Moreno for his support during the last four years.
The authors want to thank Professors  L. Fuchs, L. S. Levy and R.
Laubenbacher for their comments, support and helpful suggestions while preparing this paper.}

\begin{abstract}
We introduce a general constructive method to find a $p$-basis (and the Ulm invariants)
of a finite Abelian $p$-group $M$. This algorithm is based on Gr\"obner bases  theory.
We apply this method to determine
the additive structure of indecomposable modules over the following
Dedeking-like rings:
$\Z C_p$, where $C_p$ is the cyclic group of order a prime
$p$, and the $p-$pullback $\{\Z \rightarrow \Z_p
\leftarrow \Z \}$ of $\Z \oplus \Z$.
\end{abstract}

\maketitle
\section{Introduction}

Let $R$ be an algebra. Finding the additive structure of an $R$-module as an
Abelian group associated to a representation is a classical problem
solved in a similar way to obtaining the Jordan canonical form
of a matrix over a field, see \cite[Chapter III]{CR} and \cite[Chapter 12, \S 2]{A}.
This information is used, for example, to
determine the matrices associated to the group representation. This is
accomplished by finding a $p$-basis for the torsion part of the
group that permits a unique matrix representation for this Abelian
finite $p$-group. In 1949, Szekeres started the classification and matrix
description of modules over $\Z_{p^n} C_p$.
Since then it has been studied in detail, see
\cite{AL,AB, L1,L2,NR,NRSB,Sz}.
In \cite{L1, L2}, Levy studied these modules in the more general
context of modules over a pullback of two
Dedekind rings with a common field, which he called \emph{Dedekind-like}
rings.

Until now, the simplest way to find the additive structure of an $R$-module
consists in writing the relations as a matrix with entries in $\Z$, performing elementary
transformations over an Euclidean domain (like $\Z$), and using
the division algorithm to write the matrix in a canonical
form, see \cite[Theorem 16.8]{CR}. This approach becomes rather difficult
when the generating set is not minimal and there are several
relations among the generators. In here, we present a different
method that has the advantage of producing different group presentations
by writing the relations as polynomials and changing the term
orders used to reduce them. Furthermore, we show how to use this procedure to
find a good $p$-basis which gives the \emph{Ulm invariants} \cite{F} of
$M$ and also the \emph{type} of $M$.

The main contribution of this paper is a
constructive method to find a $p$-basis (and the Ulm invariants)
of a finite Abelian $p$-group $M$ from a given presentation of $M$
encoding the action of $p$.
The algorithm  is obtained by noting that there are some
invariant properties between the order of elements in
an Abelian group and the basis elements of certain \emph{toric ideals} \cite{S}. To accomplish this,
we use several tools from \emph{Gr\"obner bases}
\cite{Bu,CLO} and \emph{chain-modules} \cite{AB}. Furthermore, this method can be
used in general for modules over algebras on $\Z_{p^n}$ and $\Z$.

Let $M$ be a finitely generated Abelian group. We assume that $M$ is
also finitely presented, that is, $M=\la C,R\ra$, where $C$ is a
nonminimal finite generating set, and $R$ is a finite set of
relations among the elements of $C$, see \cite{G}. For example

 \[M = \la c_1, \dots , c_n \mid \sum_{j=1}^q a_{ij}c_j  = 0
\mbox{ with } a_{ij}\in \Z, \mbox{ for all } i=1,\dots,n \mbox { and } j = 1,\dots, q\ra.\]

We want to find the \emph{torsion-free rank} of $M$ and the
\emph{Ulm invariants} of the $p_j$-Sylow subgroups of $M$. This is
an old problem, the new aspects in this work are: (1) we use
Gr\"obner bases to solve the problem, and (2) using the notation and
classification introduced by Levy, we apply this method to determine
the additive structure of indecomposable modules over certain
Dedekind-like rings. In this case, the algorithm computes a
$p$-basis for the torsion part of the group.

This paper is organized as follows: in Section 2,
we introduce \emph{toric ideals} associated  to finitely
generated Abelian groups. In Section 3,
we give a description of the
\emph{reduced Gr\"obner basis} \cite{CLO} of a toric ideal
associated to a finitely generated Abelian $p$-group.
As a consequence, in Section 4, we obtain
an algorithm to compute the $p$-basis and the type of any finite
Abelian $p$-group. As an application of this algorithm, in Section 5, we show how
 to obtain the additive structure of any indecomposable module over
$\Z C_p$, where $C_p$ is the cyclic
group of order a prime $p$ and over the $p-$pullback $\{\Z \rightarrow \Z _p
\leftarrow \Z \}$ of $\Z \oplus \Z$.

\section{Gr\"obner bases associated to finitely generated Abelian groups}

We start by reviewing some concepts in finitely generated Abelian
group theory.

\begin{definition}[type]
If $p_1<\cdots <p_r$ and $M$ is a finitely generated Abelian group,
such that
\[M\cong \Z^{s_0}\oplus({\Z_{p_1}}^{s_{11}}\oplus\cdots \oplus
{\Z_{p_1^{n_1}}}^{s_{1n_1}})\oplus \cdots
\oplus({\Z_{p_r}}^{s_{r1}}\oplus\cdots \oplus
{\Z_{p_r^{n_r}}}^{s_{rn_r}}) ,\] as an Abelian group, then
the \emph{type} of $M$ is
\[\underline{t}(M)=(s_0,s_{11},\ldots, s_{1n_1},\ldots, s_{r1},\ldots,
s_{rn_r})\in \Z^{n_1+\cdots +n_r+1}.\]
\end{definition}

The number $s_0$ is the \emph{torsion-free rank} of $M$,  the number $n_i$
is the \emph{torsion rank} of $M_i$, and the sequence of numbers $s_{i1},\
\ldots,\  s_{in_i},$ are the \emph{Ulm invariants} of the $p_i$-Sylow
subgroup $M_i={\Z_{p_i}}^{s_{i1}}\oplus\cdots \oplus
{\Z_{p_i^{n_i}}}^{s_{in_i}}$ of $M$.

\begin{definition}[$p$-basis]\label{pB}
If $M$ is a $p$-group, for some prime number $p$, a set
$B=\{b_1,\ldots,b_d\}\subset M$ is called a $p$-basis of $M$ if
$M\cong_{\Z_{p^n}}\la  b_1\ra \oplus \cdots \oplus \la  b_d\ra $. A set $B$ is
a $p$-basis  of $M$ if and only if, for all $m\in M$ the sum
$m=\sum_{i=1}^d l_ib_i$ is unique,  where $0\leq l_i\leq
\ord(b_i)$ and $\ord(b_i)$ is the order of $b_i$ in the group $M$, see \cite{F}.
\end{definition}

Let $M= \Z^{s}\bigoplus_{t=1}^r M_t$ be  a finitely generated
Abelian group, where $M_t$ is the $p_t$-Sylow subgroup of $M$ with
$p_t$-rank equal to $d_t$. Consider a nonminimal  generating set
$C_t$ of each $M_t$,  such that, $C_i \cap C_j = \emptyset$ for all
$i\neq j$, and a generating set $C_0$ of $\Z^{s}$. If $C=\cup
_{t=0}^r C_t=\{c_1,\ldots ,c_q\}$, where $q\ge \sum_{t=1}^rd_t+s$,
then $\la C\ra=M $. Consider the semigroup homomorphism \[\gamma  :
\N^q \longrightarrow M, \qquad  v = (v_1,\ldots,v_q) \longmapsto
\sum_{i=1}^qv_i c_i.\]

Let $k$ be an infinite field. The previous map lifts to the following short exact sequence

\begin{equation}\label{ses3}
 0\longrightarrow \Ker (\widetilde{\gamma})\longrightarrow k[\mathbf{x}]\overset{\widetilde{\gamma}}{
 \longrightarrow} k[M] \longrightarrow 0,
 \end{equation}

 where  $k[\x]=k[x_1,\ldots,x_q]\cong k[\N^q]$  is the polynomial ring in $q$
indeterminates over $k$.  The
monomials in $k[{\bf x}]$ are denoted by ${\bf x}^{\bf
a}=x_1^{a_1}\cdots x_q^{a_q}$, where ${\bf a}=(a_1, \ldots ,a_q) \in \N^q$.
On the other hand, if $d=\sum_{t=1}^rd_t$, we have the following isomorphism
\[   k[M]\cong k[{\bf t}] = k[t_1,\ldots ,
t_d,t_{d+1},t^{-1}_{d+1},\ldots ,t_{d+s},t_{d+s}^{-1}]/\la
t_1^{k_1}-1,\ldots ,t_d^{k_d}-1\ra,\]

where $k_i$ is the order of the corresponding element in the external direct sum of $M$. Furthermore, in this external direct sum, the element $c_i\in C$ can be expressed as a tuple
\(\bfc_i=(c_{i1},\ldots,c_{id}, c_{id+1},\ldots ,c_{id+s}).\)  So,
we have a  homomorphism  of semigroup algebras

\begin{equation*}
 \widetilde{\gamma} : k[\mathbf{x}]\longrightarrow k[\mathbf{t}], \qquad x_i\longmapsto \mathbf{t}^{ \bfc_i}
=t_1^{c_{i1}}\dotsm t_d^{c_{id}}t_{d+1}^{c_{id+1}}\dotsm t_{d+s}^{c_{id+s}}.
 \end{equation*}

We denote the kernel of $\widetilde{\gamma}$ by $I_C$.
We will show how to obtain  a minimal generating set that is a $p$-basis
of $M$, from a certain Gr\"obner basis of this ideal.
In the following, we assume that we have a term order $\prec$
defined in $k[{\bf x}]$. Then every nonzero polynomial $f\in k[{\bf
x}]$ has a unique initial monomial, denoted $in_{\prec}(f)$.
Observe that for any $v=(v_1, \ldots ,
v_q)\in  \Z ^q$,  we can  write $v=v^+ - v^-$, where
$v^+=(v_1^+, \ldots ,  v_q^+)$ and $v^-=(v_1^-, \ldots , v_q^-)$ are
nonnegative integer tuples. Denote by
 $\Ker(\gamma  )$  the subgroup of $\Z^q$ consisting of all
elements  $ v$ such that ${\gamma }(v^+)={\gamma } (v^-)$. Let

\[P(\Ker (\gamma))=\big\{{\bf x}^{ v^+}- {\bf x}^{ v^-} \mid v\in \Ker(\gamma
)\big\}.\]

The following lemma follows immediately from \cite[Lemma 4.1]{S}.

\begin{lemma}
The  ideal $I_C$ is generated as a $k$-vector space by the set
$P(\Ker (\gamma))$.
\end{lemma}

 Recall that $C$ is a nonminimal finite generating set  of $M$.
We assume that there exists a finite set of defining relations $\R$ for $C$ in $M$.
We use the notation $ v^+$ and $ v^-$ to
write the relations as
 $\sum_{t=1}^qv^+_tc_t=\sum_{t=1}^qv^-_tc_t$. These relations induce a subset of vectors in $\Z^q$ and
 a subset of polynomials  in $I_C$

  \begin{align*}
  \overline\R &=\big\{ v\in
\Z^q: \sum_{t=1}^qv^+_tc_t=\sum_{t=1}^qv^-_tc_t \hbox { is in
}\R\big\}\subset \Z^q,\\
P(\overline\R ) &=\big\{ P(v)={\bf x}^{v^+}-{\bf
x}^{v^-}\mid v\in\overline \R\big\}\subset I_C.
\end{align*}

Let $G_{\overline \R}$ denote the reduced Gr\"obner basis of
the ideal  generated by $P(\overline\R)$,
with respect to  the order $\prec$. Similarly, this Gr\"obner basis induces the  set
 $\overline\R_\G$ of tuples

\[\overline\R_\G=\big\{v\in \Z^q \mid {\bf x}^{v^+}-{\bf x}^{v^- }\in
\G_{\overline \R}\big\}\]
and the set $\R_\G$ of relations
\[R_{\G} = \big\{ \sum_{t=1}^qv^+_tc_t=\sum_{t=1}^qv^-_tc_t \mbox{ such that } v\in
\overline\R_\G\big\}.\]

\begin{proposition} \label{I}
Let    $ v_k\in\overline \R$ with $in_{\prec}(P( v_k))={\bf x}^{
v_k^+}$ for $k\in{1,2}$. Also let $w_1,w_2\in \N^q$ such that
  $w_1+v_1^+=w_2+v_2^+$.
 If $P(v)={\bf x}^{ w_1}P(v_1)-{\bf x}^{w_2}P(v_2)$ then
$v=v_1-v_2$.
\end{proposition}

\begin{proof}
  We have  $P( v)={\bf x}^{ w_1}P(v_1)-{\bf
x}^{w_2}P(v_2)= {\bf x}^{ w_2+ v^-_2}-{\bf x}^{ w _1+ v^-_1}.$ But
\[\gamma ( w_1+ v^+_1)= \gamma (w_2+ v^+_2)=\gamma ( w_1+
v^-_1)=\gamma ( w_2+ v^-_2),\] \[ and \  v=( w_2+ v^-_2)-( w_1+
v^-_1)=( w_2+ v^-_2)-( w_1+ v^-_1)+( w_1+
v^+_1)-(w_2+v_2^+)=v_1-v_2.\]
\end{proof}

 \begin{theorem}\label{th:ge}
The set $\R_\G$ is  a set of relations  for $C$ in $M$.
 \end{theorem}

 \begin{proof}
First observe that,
if $v\in \Ker(\gamma )$, then $\sum_{t=1}^qv_tc_t=0$. So  $v\in \la
\overline \R \ra$
 and since $\la \overline \R \ra \subset \Ker(\gamma )$, then
the set $\overline \R $ generates the subgroup $\Ker(\gamma )$.
This implies that
the ideal $I_C$ is generated by $P(\overline\R )$.
Next, we use the
Buchberger algorithm \cite{Bu} to obtain the  reduced Gr\"obner
basis of $I_C$ from $P(\overline\R)$.

 By Proposition \ref{I}, we have that the $S$-polynomial $S(P(u_1),P( u_2))=P(v)$ satisfies $ v = u_1-u_2\in \Ker (\gamma).$ Let
$S$ be the set of all nonzero $S$-polynomials obtained in the
Buchberger algorithm and let $\overline S =\{ v\in \Z ^q \mid P(
v)\in S \}.$ Clearly $\la \overline\R \ra=\la \overline\R \cup \overline S \ra$.
We denote by $\overline{\R'} =\overline\R \cup \overline S$.
Now, we reduce the set of polynomials in $P(\overline{\R'}) $. Suppose
$in_{\prec}(P( v))={\bf x}^{ v^+}$  divides $in_{\prec}(P(
v_1))={\bf x}^{ v_1^+}$, with $ v,v_1\in  \overline{\R'}$.
There exists $ w\in \N^q$ such that $w+v^+=v_1^+$  and
$P(v_2)=P(v_1)- {\bf x}^{w}P(v)$. Hence $v_2=v_1-v\in
\overline{\R'}$ by Proposition \ref{I}. Then $\la \overline\R
\ra=\la \overline{\R'} \ra=\la (\overline{\R'}\setminus\{ v_1\})\cup
\{v_2\}\ra$. So $\la \overline\R_\G \ra=\la \overline\R \ra$. This
proves our claim.
\end{proof}

\section{The reduced presentation of $M$}

Given a generating set $C$ of a finite Abelian $p$-group $M$, one can obtain a
set of relations by studying the action of $p$ over the elements in $C$.
In the last section, we saw that any Gr\"obner basis of $I_C$ gives a set of relations
for $M$. In this section, we describe a particular Gr\"obner basis that gives a
$p$-basis of $M$.
 We assume that the elements of $C$ have orders
$\ord(c_1)\geq \cdots \geq \ord (c_q)$. Consider the following  chain of
subgroups of $M$
\[\la c_1\ra \subseteq \la c_1,c_2\ra \subseteq \cdots \subseteq \la c_1,c_2,\ldots,
c_s\ra \cdots\subseteq \la C \ra = M.\]

For $s \ge 2$, let $r_{s}=\min\{k \mid p^k
c_s\in \la c_1,c_2,\ldots ,c_{s-1}\ra\}$.
There are two possibilities, either $r_{s} < \ord(c_s)$, or  $r_{s} = \ord(c_s)$, in this case,
$p^{r_{s}}c_s=0\in \la c_1,c_2,\ldots ,c_{s-1}\ra$. Thus, we have the following
set of relations

\[R_p= \big\{\ p^{r_1}c_1=0,\ p^{r_s}c_s=p^{r_s}\sum_{t<s}a_{st}c_t,\ where  \ a_{st}\in
\Z,\ for \ 2\leq s\leq q\big\}.\]

\begin{proposition}\label{pro:R} The  relations $R_p$
together with the set $C$ is a presentation of the $p$-group $M$.
\end{proposition}

\begin{proof} Suppose $\sum_{t=1}^q \ell _t c_t=0.$ Dividing  $\ell _q $ by
$p^{r_q}$, we obtain  $\ell _q=s_qp^{r_q}+s'_q$ with $0\le
s'_q<p^{r_q}$. Therefore,  $s_qp^{r_q}c_q+s'_qc_q+\sum_{t=1}^{q-1}
\ell _t c_t=0$.
Suppose that $s'_q\ne 0$. Using the relation
$p^{r_q}c_q=\sum_{t<q}p^{r_q}u_{qt}c_t$, we obtain $\sum_{t=1}^q
\ell _t c_t=s'_qc_q+\sum_{t=1}^{q-1} \ell '_t c_t=0$. If  $\gcd
(s'_q,p)=1$, then  $c_q=\sum_{t<q-1}l_tc_t$ and $\sum_{t=1}^q \ell
_t c_t=\sum_{t=1}^{q-1} l_t c_t=0$.
 If $\gcd (s'_q,p)=p$, let $p^r$
be the maximum number such that $p^r$ divides $s'_q$. Then
$p^rc_q=-\sum_{t=1}^{q-1} \ell '_t c_t$. But this  is impossible,
because $p^{r_q}$ is the minimum with this condition. Thus
$\gcd (s'_q,p)=1$.

Repetition of this argument shows that $\sum_{t=1}^q \ell _t
c_t=l'_1 c_1=0$ with $\gcd (l'_1,p)=1$. But this implies that
$c_1=\ldots =c_q=0$ which  is impossible. Then $s'_q=0$ and
the relation $\sum_{t=1}^q \ell _t c_t=0$ is a linear combination
of the relations in  $\R$.
\end{proof}

\begin{proposition}\label{prop:rGb}
Let $\prec$ be the lexicographic ordering with $x_1\prec x_2 \prec\cdots \prec x_q$.
Then,
the reduced Gr\"obner basis of $I_C$ with respect to $\prec$ equals
\[{\G}_p=\big\{x_1^{p^{r_1}}-1, x_2^{p^{r_2}} - x_1^{a_{21}p^{r_2}},\ldots ,x_q^{p^{r_q}}-
\prod_{t=1}^{q-1}x_t^{a_{qt}p^{r_q}}\big\}.\]
\end{proposition}

\begin{proof}
Observe that $\G_p = P(\overline{\R_p})$. Thus,  by Theorem \ref{th:ge},
$\G_p$ generates $I_C$. Furthermore, $\G_p$ is a reduced Gr\"obner basis
since $\gcd(in_{\prec}(p_1), in_{\prec}(p_2)) = 1$, for any $p_1$ and $p_2$ in $\G_p$. This forces all S-polynomials to
be zero modulo $\G_p$, see \cite{CLO}.
\end{proof}

Let $pR$ be  the following set of relations for $C$ in $M$

\begin{equation}\label{pR}
pR= \big\{ pc_q=0, \quad pc_j=\sum_{t>j} a_{jt}c_t, \mbox{ for all } 1\leq j\leq q-1,
\mbox{ with } 0\leq a_{jt} \in \Z \big\}.
\end{equation}

 These relations can be used to find  the order of any element in $M$,
 since $pM$ is the \emph{Frattini} subgroup of $M$, \cite{KM} . On the other hand,
 from the action of $p$, we can find the minimal number of
generators of $M$, that is, the $p$-rank of $M$ by
 the Burnside Basis Theorem for finite groups ($M/pM$), \cite{H}.
Let $d$ be the $p$-rank of $M$. For each $t\geq 2$, let $D_t$ be the
set
\[D_t = \big\{c_t-\sum_{j=1}^{t-1}a_{ts}c_j \mid 0\leq a_{tj}\leq\ord (c_j)\big\}.\]

If $b_t$ is the element of maximal order in $D_t$, then
$p^{r_t}=\ord (b_t)$.  Therefore, we have the following set of relations, denoted by
$R_{pbasis}$

\[ \big\{ p^{r_1}c_1=0, \
p^{r_t}c_t=p^{r_t}\sum_{j<t}a_{tj}c_s \mbox{ for } 2\leq t\leq d,  \mbox{ and }
c_t=\sum_{j<d}a_{tj}c_j \mbox{ for } d< t\leq q \big\}.\]

It is clear that $M = \la b_1,b_2,\ldots ,b_d\ra$, so
$R_{pbasis}$ is a set of relations for $C$ in $M$. As a corollary of Proposition
\ref{prop:rGb},  we have

\begin{corollary}\label{cor:pG}
Let $\prec$ be the lexicographic ordering with $x_1\prec x_2 \prec\cdots \prec x_q$.
Then, the reduced Gr\"obner basis of $I_C$ with respect to $\prec$ equals

\begin{equation}\label{Gpbasis}
\G_{pbasis} = \big\{x_1^{p^{r_1}}-1,\ldots
,x_d^{p^{r_d}}-\prod_{t=1}^{d-1}x_t^{a_{dt}p^{r_d}},
x_{d+1}-\prod_{t=1}^{d}x_t^{a_{td+1}},\ldots,x_q-\prod_{t=1}^{d}x_t^{a_{tq}}\big\}.
 \end{equation}
\end{corollary}

Note that $\G_{pbasis}$ is just a refinement of $\G_p$ obtained by setting some of the
$p^{r_j}$ equal to 1. The next theorem is the key to our algorithm. It says that
the generating set obtained from $\G_{pbasis}$ is actually a $p$-basis of $M$.

\begin{theorem}\label{th:pre}
The set $\overline{C}=\{c_1,\ldots,c_d-\sum_{t=1}^{d-1}a_{dt}c_t\}$
is a $p$-basis of $M$.
\end{theorem}

\begin{proof}
We have seen that  $M=\la b_1,b_2,\ldots
,b_d\ra = \la \overline{C} \ra$. Now, we  will prove that the sum $\la b_1\ra + \dotsm+ \la
b_d\ra$ is actually a direct sum.  If $y\in \la b_t\ra \cap \la
b_1,\ldots,b_{t-1}\ra$. Then $y=\alpha _tb_t=\sum_{j=1}^{t-1}\alpha
_sb_s$. Thus, $\alpha _tc_t=\sum_{j=1}^{t-1}\alpha' _jc_j$.
The argument preceding this theorem shows that
$\alpha_t\ge p^{r_t}$, so $\alpha _t=\alpha'_tp^{r_t}+\beta _t$,
with $0\le \beta _t<p^{r_t}$. Thus $\beta
_tc_t=\sum_{j=1}^{t-1}\alpha " _jc_j$ which implies $\beta _t=0$ and
$y=\alpha '_tp^{r_t}b_t=0$. So $M=\bigoplus_{t=1}^d \la b_t\ra$.
\end{proof}

We can summarize the above results as follows.

\begin{remark} \label{rmk:1}\mbox{}
\begin{itemize}
\item[(1)] Given a presentation of a finite Abelian $p$-group $M = \la C, \R\ra$, there
exists a term ordering such  that the reduced Gr\"obner basis of the toric ideal $I_C$ gives a $p$-basis for $M$.
\item[(2)] Given a homomorphism $\widetilde{\gamma}$ as in \eqref{ses3}.
The presentations for the corresponding finite Abelian group $M$  can be obtained
from Gr\"obner bases of the toric ideal $\Ker(\widetilde{\gamma})$.
\end{itemize}
\end{remark}

\section{The $p$-basis algorithm and the additive structure of $M$}

Corollary \ref{cor:pG} gives an explicit description of a reduced
Gr\"obner basis for $I_C$. Moreover, Theorem \ref{th:pre} shows that
the corresponding set of generators is a $p$-basis of $M$.
Nevertheless, we obtained this Gr\"obner basis from a very special
set of  relations whose definition was nonconstructive, namely
$\R_{pbasis}$.   In particular, this set of relations specified the
ordering on the indeterminates for the specific lexicographic order
needed in Corollary \ref{cor:pG}. In this section, we put all these
results together to compute the invariants of a finite Abelian
$p$-group $M$ from a particular presentation.

Let $pR$ be the finite presentation of $M$ introduced in \eqref{pR}, that is, assume
that the action of $p$ in a generating set $C$ is known. Following Remark \ref{rmk:1},
we need to find an ordering of the indeterminates, such that, the Gr\"obner basis with
respect to the corresponding lexicographic order has the form \eqref{Gpbasis}.
Note that there might be several such orderings. In the last section, we saw that
if $\ord(c_i) < \ord(c_j)$ then $x_j \prec x_i$. We also need to break ties among the elements in $C$ with the same order in the group.

In practice, one first break ties arbitrarily. If the Gr\"obner basis has the required
form, we are done. Otherwise, there is an element in the Gr\"obner basis of the form
$x_j^{p^{r_j}} - x_i \prod_{t=1, t\neq i}^{j-1}x_t^{a_{jt}p^{r_j}}$, with
$\ord(c_j) = \ord(c_i)$. In this case, we need to invert the order of $x_i$ and $x_j$ to
$x_j < x_i$. This process eventually terminates, moreover; it effectively gives the
desired Gr\"obner basis since the $p$-basis itself always exists.
The output of the algorithm
consists on the $p$-basis and the Ulm invariants of $M$, that is, the type $\underline
t(M)$.

\begin{algorithm}\label{alg:pb}
    Input: $C$, $pR$.

\begin{itemize}
\item[(A1)] Write the relations in $pR$ as polynomials in $k[\x]$
as follows:  $x_q^p -1,$ and $x_j^p-\prod_{t>j}x_t^{a_{jt}}$, for
$1\leq j\leq q$.

\item[(A2)] Find the order of all $c_j$ by computing all the univariate polynomials
in the ideal $I$ generated by the polynomials obtained in (A1).

\item[(A3)] Find an ordering of the indeterminates, such that, the reduced
lexicographic Gr\"obner basis $\G_p$ of $I$ has the form \eqref{Gpbasis}.

\item[(A4)] Let $d$ be the number of polynomials in $\G_p$ such that the
initial term has exponent $>1$.
If $p^{r_j}>1$ and $x_j^{p^{r_j}}-\prod_{t=1}^{j-1}
x_t^{a_{jt}p^{r_j}}\in \G_p$,  then add $b_j$ to the
$p$-basis, where $b_{j}$ is the following element of order $p^{r_j}$:
\[b_j=c_j-\sum_{t=1}^{j-1}a_{jt}c_t. \]

\item[(A5)] To compute the type of $M$, let $s_{r}$ be
the number of elements with the same order $p^r$.
Then $\underline{t}(M)=(s_1,\ldots,s_n)$.
\end{itemize}
Output: $B=\{b_1,\ldots ,b_d\}$ and $M\cong (\Z_p)^{s_1} \oplus
\dotsm \oplus (\Z_{p^n})^{s_n}$.
\end{algorithm}

\begin{example}
Let $M=\la c_1,c_2,c_3,c_4,c_5,c_6,c_7,c_8\ra$
be a $5$-group, with the following relations:
\begin{align*}
&5c_1-c_8-4c_5-2c_6-3c_7=0, \quad 5c_2-4c_6-2c_7=0, \quad 5c_3-4c_7=0,\\
& 5c_4=0, \quad 5c_5=0, \quad 5c_6=0, \quad 5c_7=0, \quad 5c_8=0.
\end{align*}

The corresponding polynomials are

\[\big\{x_1^5-x_8x_5^4x_6^2x_7^3,x_2^5-x_6^4x_7^2,x_3^5-x_7^4,x_4^5-1,x_5^5-1,x_6^5-1,x_7^5-1,x_8^5-1\big\}.\]

The reduced lexicographic Gr\"obner basis equals

\begin{align*}
\big\{
&x_1^{25}-1, \quad x_2^{25}-1, \quad x_3^{25}-1, \quad x_4^5-1, \quad x_5^5-1 \\
&x_6-x_3^{15}x_2^{20}, \quad x_7-x_3^{20}, \quad x_8-x_5x_3^{10}x_2^{10}x_1^5
\big\}.
\end{align*}

In this case, $d=5$. So, the Gr\"obner basis gives the following information

\begin{align*}
&25c_1=0, \quad 25c_2=0, \quad 25c_3=0, \quad  5c_4=0, \quad 5c_5=0, \\
&c_6=15c_3+20c_2, \quad c_7=20c_3, \quad c_8=c_{5}+10c_3+10c_2+5c_1.
\end{align*}

Hence, the  $p$-basis is equal to
$B=\{c_1,c_2,c_3,c_4,c_5\}$, $M\cong {\Z_5}^2\oplus {\Z_{25}}^3$,  and
$\underline{t}(M)=(2,3)$.
\end{example}

The classical way to solve this
problem, using matrix transformations over an Euclidean domain,
appears in \cite{A}. We remark that it is possible to perform the second
step in the algorithm because by definition, $I_C$ is a zero-dimensional ideal.
Moreover, each univariate polynomial in $I_C$ has the form $x_j^{\ord(c_j)}-1$.

\section{Indecomposable modules over Dedekind-like rings}

 Let $R_1$ and $R_{2}$ be two rings. Let $R$ be the
pullback ring of the rings $R_i$ over a common ring $\overline R$,
that is, $R=\{R_1\rightarrow \overline R \leftarrow R_2\}$. In \cite
{L1}, L. Levy studied the separated representation of an $R$-module
$M$. In \cite {L2}, he described the indecomposable $R$-modules
when $R_1$ and $R_2$ are Dedekind domains and $\overline R $ is a
field $k$ ($R$ is  called a Dedekind-like ring). In
particular, he studied modules over two rings: $\Z C_p$, where $C_p$
is the cyclic group of order a prime number $p$, and the
$p-$pullback $\{\Z \rightarrow \Z _p \leftarrow \Z \}$ of $\Z \oplus
\Z$.

An $R$-module $S$ is separated if it is an $R$-submodule of a
direct sum $S_1\oplus S_2$, where each $S_i$ is an $R_i$-module. A
separated representation of an $R$-module $M$ is an $R$-module
epimorphism $\phi :S\longrightarrow M$, such that,  $S$ is a separated
$R$-module and  if $\phi$ admits a factorization $\phi
:S{\overset {f}{\longrightarrow}}S'\longrightarrow M$ with $S'$ also a
separated $R$-module, then   $f$ must be one to one. Let $P_i=\ker
(R_i\longrightarrow k)$, then $P=\{P_1 \rightarrow 0 \leftarrow P_2\}$
is an ideal of $R$ . We call an $R$-module $M$ $P$-mixed, if each
\emph{torsion} element $m$ is annihilated by some power of $P$.
The separated modules $S=\{S_1
{\overset{f_1}{\longrightarrow}}k{\overset{f_2}{\longleftarrow}}S_2\}$
satisfying one of the following two conditions: (1) $ S_i\cong \hbox{ nonzero ideal of }
R_i$, or  (2) $ S_i \cong R_i/P^{e}_i$
 form the  basic building
blocks for all finitely generated, $P$-mixed  $R$-modules. If $S$
is a building block, then $S$ has exactly one submodule which has
the form $\{X\rightarrow 0\leftarrow 0\}$ and is $R$-isomorphic to
$k$ (left $k$ of $S$). Similarly, $S$ has a
 right $k$ of $S$.

\begin{definition} (Deleted Cycle and Block Cycle Indecomposables)\\
\begin{itemize}
\item[(a)] Let $S^{(1)},\ldots,S^{(m)}$ be a sequence of basic building
blocks, such that,
\[\begin{CD}
S^{(1)} @>>> S_{21} \\
@VVV            @V\sigma_{21}VV \\
S_{11} @>\sigma_{11}>> k
\end{CD}
\dotsm
\begin{CD}
S^{(i)} @>>> S_{2i} \\
@VVV            @V\sigma_{2i}VV \\
S_{1i} @>\sigma_{1i}>> k
\end{CD}
\dotsm
\begin{CD}
S^{(m)} @>>> S_{2m} \\
@VVV            @V\sigma_{2m}VV \\
S_{1m} @>\sigma_{1m}>> k
\end{CD}
\]
and suppose that for $1\leq i\leq m$, $S^{(i)}$ has a right $k$ and
$S^{(i+1)}$ has a left $k$.  A deleted cycle indecomposable $M$
is the direct sum $S=\bigoplus _{i=1}^m S^{(i)}$ modulo a relation
which identifies the right $k$ of $S^{(i)}$ with the left $k$ of
$S^{(i+1)}$, that is, first choose $p_j\in P_j-P_j^2$ for $j=1,2$,
then make the following identification
\[p_2^{d(2,i)-1}s_{2,i}=-p_1^{d(1,i+1)-1}s_{1,i+1},\mbox{ where }  s_{j,i}\in
S_{ji},\ \sigma _{ji}(s_{j,i})= \overline 1,\]
for $j=1,2$; $1\leq i\leq m-1$, and $d(j,i)$ the length of $S_{ji}$. In other words,
it is the direct sum $S$ modulo
\[\big\{p_2^{d(2,1)-1}s_{2,1}+p_1^{d(1,2)-1}s_{1,2}, \ \ldots , \
p_2^{d(2,m-1)-1}s_{2,m-1}+p_1^{d(1,m)-1}s_{1,m}\big\}.\]

\item[(b)] Let $S^{(1)},\ldots,S^{(m)}$ be a sequence of basic building
blocks

\[\begin{CD}
S^{(1)} @>>> S_{21} \\
@VVV            @V\sigma_{21}VV \\
S_{11} @>\sigma_{11}>> k
\end{CD}
\dotsm
\begin{CD}
S^{(i)} @>>> S_{2i} \\
@VVV            @V\sigma_{2i}VV \\
S_{1i} @>\sigma_{1i}>> k
\end{CD}
\dotsm
\begin{CD}
S^{(m)} @>>> S_{2m} \\
@VVV            @V\sigma_{2m}VV \\
S_{1m} @>\sigma_{1m}>> k
\end{CD}
\]

 each with a left and a right $k$. Write $m=l\underline m$, where
$\underline m$ is the unique smallest positive integer, such that,
for all $i$, $S^{(i)}\cong S^{(i+\underline m)}$. Let
$f(z)=\lambda_o+\lambda_1z+\dotsm + \lambda_{l-1}z^{l-1}+z^l$ be a
power  of an irreducible polynomial  in $k[z]$. A block cycle
indecomposable $M$ is a  deleted cycle indecomposable  modulo  the
following relation $
-p_2^{d(2,m)-1}s_{2,m}=\sum_{j=0}^{l-1}\lambda_jp_1^{d(1,j)-1}s_{1,(j\underline
m+1)},$ which identifies the right $k$ of $S^{(m)}$ with a
one-dimensional subspace of $S_{1,1}\oplus S_{1,\underline m
+1}\oplus S_{1,(2\underline m+1)}\oplus \dotsm$. In other words, it is the direct
sum  $S$ modulo
\[
\begin{split}
\big\{p_2^{d(2,1)-1}s_{2,1}+p_1^{d(1,2)-1}s_{1,2}, \ \ldots, \
&p_2^{d(2,m-1)-1}s_{2,m-1}+p_1^{d(1,m)-1}s_{1,m}, \ \\
&p_2^{d(2,m)-1}s_{2,m}+\sum_{j=0}^{l-1}\lambda_jp_1^{d(1,j)-1}s_{1,(j\underline
m+1)}\big\}.
\end{split}\]
\end{itemize}
\end{definition}

As a consequence, if $M$ is a
deleted cycle then  $1\leq d(2,i)\ne \infty$, for $1\leq
i\leq m-1$ and $1\leq d(1,i)\ne \infty$, for $2\leq i\leq m$. But
the length of either one of  $S_{11}$ or $S_{2n}$ may be infinite.
If $M$ is a block cycle, then
$1\leq d(j,i)\ne \infty$ for $1\leq i\leq m$ and $j=1,2$.

\begin{remark}
The indecomposable, finitely generated, $P$-mixed modules are
deleted cycle indecomposables and block cycle indecomposables.
Every separated $R$-module is a direct sum of basic building blocks.
Moreover, basic building blocks are always indecomposable
$R$-modules, see \cite{L2}.
\end{remark}

\subsection{Additive descriptions}

Using Algorithm \ref{alg:pb}, we describe the
additive structure of the indecomposable $R$-modules when $R$ is one
of the following rings:  $\Z C_p$ or the $p$-pullback of $\Z \oplus
\Z$,
 $\{ \Z \rightarrow\Z_p \leftarrow \Z\}$.
 In these two  cases the concept of $P$-mixed coincides with
$p$-mixed.

{\bf The ring $\Z C_p$}:

\[ \begin{CD}
\Z C_p @>>> \Z[\zeta ] \\
@VVV            @V\nu_{2}VV \\
\Z @>\nu_{1}>> \Z_p
\end{CD}\]

 Let $\zeta$ be a primitive $p$th root of
unity, and let $x$ be a generator of $C_p$. Then $\Z C_p\cong \{ \Z
{\overset{\nu _1 }{ \longrightarrow}}\Z_p {\overset {\nu _2 }{
\longleftarrow}}\Z [\zeta ]\},$ where the isomorphism is given by $
x\longrightarrow (1\rightarrow \overline 1\leftarrow \zeta).$ The action
of $p_1$ and $p_2$ in $\Lambda=\Z C_p$ is given by the following
formulas:

\begin{align*}
p_1 &= x^{p-1}+x^{p-2}+\cdots + x+1=\big\{p\rightarrow 0\leftarrow
\zeta ^{p-1}+\zeta ^{p-2}+\cdots +\zeta +1\big\} \mbox{  and } \\
p_2 &= x-1=\{0\rightarrow 0\leftarrow \zeta -1\}, \quad  p=(p,p)=
p_1+p_2^{p-1}\sigma (p_2), \quad p_1p_2=0.
\end{align*}

where $\sigma (p_2)$ is a
polynomial in $p_2$, with degree less or equal than $p-1$, which exists
because the sum equals $p$.  So, every
element $m$ of a $\Z C_p$-module $M=\la a_1,\ldots , a_n\ra$ is a
linear combination of these generators and the elements resulting
from the action of $p_1$ and $p_2$ over them.

\begin{example}
Let $\Lambda=\Z C_3$  and  $M=\la a\ra_{\Z C_3}$ be a
deleted cycle indecomposable with $d(1)=d(2)=3,$ and
$3=p_1+2{p_2}^2$. We need to compute
the action of $p$ in $\Lambda$ over the generator $a$
to obtain a generating set for $M$ over $\Z$. This is the classical way
to begin this problem in Abelian group theory. Thus,

\[3a=p_1a+2p_2^2a,\; 3p_1a=p_1^2a,\; 3p_1^2a=0, \;
3p_2a=p_2^2a, \; 3p_2^2a=0.\]

The generating set is $C=\{a,p_2a,p_1a,p_2^2a,p_1^2a\}$, the corresponding
ideal is generated by the binomials
\[ \big\{ x_1^3-x_3x_4^2,\; x_3^3-x_5,\; x_5^3-1,\; x_2^3-x_4,\; x_4^3-1 \big\}.\]
The order of each element in $C$ is $\{27,9,9,3,3\}$.
The reduced Gr\"obner basis is equal to
\[\big\{ x_1^{27}-1,\ x_2^9-1, \ x_3-x_2^3x_1^3, \ x_4-x_2^3, \ x_5-x_1^9\big\}.\]
So, Algorithm \ref{alg:pb} outputs
\[B=\big\{a,\ p_2a\big\},\ M\cong \Z_9 \oplus \Z_{27},\ \underline t(M)=(0,0,1,1).\]
The extra zero in the type means that the torsion-free rank equals $0$.
Therefore, the action of $\Lambda$ does not change if we consider $M$ as a
module over $\Z_{27}C_3$.
\end{example}

\begin{example}
Let $\Lambda=\Z C_3$ and $M=\la a\ra_{\Z C_3}$ be a
block cycle indecomposable with $d(1)=4,\ d(2)=4,$ and $f(z)=z-2$.
The action of $p=3$ is given by $3=p_1+2{p_2}^2$. We also have the
relation  $p_1^3a=2p_2^3a$. The action of $p=3$ over $a$ is given by

\[3a=p_1a+2p_2^2a,\ 3p_1a=p_1^2a,\  3p_1^2a=p_1^3a,\ 3p_1^3a=0,
3p_2a=2p_2^3a, \ 3p_2^2a=0.\]

The generating set is
$C=\{a,p_1a,p_2a,p_1^2a,p_2^2a,p_1^3a,p_2^3a\}$. In this case, the corresponding
toric ideal is generated by
\[
\big\{x_1^3-x_2x_5^2,\; x_2^3-x_4,\; x_4^3-x_6,\; x_6^3-1,\;
x_3^3-x_7^2,\; x_5^3-1,\; x_7^3-1,\; x_6-x_7^2\big\}.
\]
The order of each element in $C$ is $\{81,27,9,9,3,3,3\}$. The reduced Gr\"obner basis
is equal to
\[\big\{
x_1^{81}-1,\ x_2^3-x_1^9,\ x_3^3-x_1^{27}, \ x_4-x_1^9,\ x_5-x_2x_1^{78}, \
x_6-x_1^{27},\ x_7-x_1^{54}\big\}. \]

Hence, Algorithm \ref{alg:pb} outputs
\[B=\big\{a,\ p_1a-3a,\ p_2a-9a\big\},\ M\cong  {\Z_{3}}^2\oplus \Z_{81},\ \underline t(M)=(0,2,0,0,1).
\]
\end{example}

{\bf The $p$-pullback ring of $\Z \oplus \Z$}:
 The  $p$-pullback of $\Z \oplus \Z$ is the
subring $\Lambda =\{ \Z \rightarrow\Z_p \leftarrow \Z\}$  of $\Z
\oplus \Z$. In this case, let $p_1=(p,0)$ and $p_2=(0,p)$. Then
$p=(p,p)=p_1+p_2$.

\begin{example}
Consider the pullback ring $\Lambda =\{ \Z \rightarrow\Z_3
\leftarrow \Z\}$ and a deleted cycle indecomposable module $M=\la
a_1,a_2\ra_\Lambda$, with $d(1,1)=3$,  $d(1,2)=3$, $d(2,1)=3$,
$d(2,2)=3$, and $-4p_2^2a_1=p_1^2a_2$. Note that the order of these elements
is $3$, since they are in the \emph{socle} $M[p]$ of $M$;
thus, the last relation is $2p_2^2a_1=p_1^2a_2$. Also $p=p_1+p_2$.
Therefore, the generators are
\[
C=\big\{a_1,\;a_2,\;p_1a_1,\;p_2a_1,\;p_1a_2,\; p_2a_2,\; p_1^2a_1,\;p_2^2a_1,\;p_1^2a_2,\; p_2^2a_2\big\}.\]
Besides the previous relation $p_1^2a_2= 2p_2^2a_1$, the relations obtained from the action of $p$ are
\begin{align*}
 &3a_1=p_1a_1+p_2a_1,\ 3p_1a_1=p_1^2a_1,\ 3p_2a_1=p_2^2a_1,\ 3p_1^2a_1=0,\
 3p_2^2a_1=0,\\
 &3a_2=p_1a_2+p_2a_2,\ 3p_1a_2=p_1^2a_2,\ 3p_2a_2=p_2^2a_2,\ 3p_1^2a_2=0,\ 3p_2^2a_2=0.
 \end{align*}
The toric ideal is generated by
\[
\begin{split}
\big\{
&x_1^3-x_3x_4,\; x_3^3-x_7,\; x_4^3-x_8,\; x_7^3-1,\; x_8^3-1,\\
&x_2^3-x_5x_6,\; x_5^3-x_9,\; x_6^3-x_{10}, \;  x_9^3-1,\; x_{10}^3-1, \;
x_9-x_8^2
\}.
\end{split}
\]

The order of each element in $C$ is $\{27,27,9,9,9,9,3,3,3,3\}$. The Gr\"obner basis is equal to
\[\begin{split}
\big\{
x_1^{27}-1,\ x_2^{27}-1,\ &x_3^9-1,\ x_4-x_3^8x_1^3,\ x_5^3-x_3^3x_1^{18},\
x_6-x_5^2x_3^6x_2^3x_1^9,\\ &x_7-x_3^3,\ x_8-x_3^6x_1^9,\ x_9-x_3^3x_1^{18},\
x_{10}-x_3^6x_2^9x_1^9
\big\}.
\end{split}
\]
Using the algorithm, we obtain the $p$-basis
\[B=\big\{a_1,\ a_2,\ p_1a_1,\ p_1a_2-p_{1}a_{1}-6a_{1}\big\}, \  M\cong \Z_{3}\oplus \Z_9 \oplus {\Z_{27}}^2,
\ \underline t(M)=(0,1,1,2). \]
\end{example}

Let $M$ be an indecomposable $R$-module and let
$S=\bigoplus _{i=1}^m S^{(i)} $ be the separated representation of
$M$. If $m=1$, then $S=\{S_{1}\rightarrow k\leftarrow S_{2}\}=\la
a\ra $ is a basic building block, and  $\length(S_j)=d(j) \neq \
\infty$, for $j=1,2$, because $M$ is $\Z_{p^n}$-free. Thus, the
subset $A=\{a, p_1a,\ldots , p_1^{d(1)-1}a,p_2a,\ldots
,p_2^{d(2)-1}a\}$ generates $S$ as an Abelian group over $\Z$.
The next theorem shows how to use Algorithm
\ref{alg:pb} to obtain the type and a $p$-basis of any
basic building block with torsion part.

\begin{theorem}\label{th:bb}
Let $S=\{S_1\rightarrow k\leftarrow S_2\}=\la a\ra$ be a basic building block. Then

\begin{itemize}
\item[(i)] If $d(1)\cdot d(2) < \infty$, then Algorithm \ref{alg:pb}
gives a $p$-basis for $S$ using  the presentation $S=\la A, pA\ra$.

\item[(ii)] If $d(j)=\infty$ for exactly one $j$, one can
obtain a basis for $S$, by adding to the input of
Algorithm \ref{alg:pb} the number $\ex
(t(S))+2$, where $\ex (t(S))$ denotes the exponent of the torsion
subgroup of $(S,+)$.

\item[(iii)] If both lengths are infinite then the Abelian group $(S,+)$ is
torsion free. In this case, the rank is
$p$, and $\{a,p_2a,\ldots , p_2^{p-1}a\}$ is a $p$-basis for $R=\Z
C_p$. If $R =\{ \Z \rightarrow\Z_p \leftarrow \Z\}$ is the
$p$-pullback of $\Z \oplus \Z$, then  either
$\{ a, p_1a\}$ or $\{ a, p_2a\}$ is a $p$-basis for $S$.
\end{itemize}
\end{theorem}

\begin{proof}
In part (iii), the case  $R=\Z C_p$ is a direct consequence of
\cite[Application 1.10]{L2} and the case $R =\{ \Z \rightarrow\Z_p \leftarrow \Z\}$
is trivial.
If $d(1)< \infty$ and $d(2)< \infty$, then $S$ is an
$R_p$-module. So, $\la A, pA \ra $ is a presentation of $S$. Hence,
applying Algorithm \ref{alg:pb}, we obtain a $p$-basis.
If $d(1) = \infty$ or $d(2) = \infty$, we change the infinite length for
 $\exp(t(S))+2$. After this, we can apply Algorithm \ref{alg:pb}
to get a $p$-basis. Using the proof
of Theorem 11.6 in \cite{L2}, we can recover  the basis for $S$. If
$R=\Z C_p$ and $d(2) = \infty$, there are $p-1$ elements of
order $\exp(t(S))+2$ in the basis, by  \cite[Application 1.10]{L2}.
These elements have infinite order and the remaining elements in
the basis form the $p$-basis for the torsion part. If $d(1)=
\infty$, then there is one element with infinite order in the
basis. If $R$ is the $p$-pullback of $\Z \oplus \Z$, we have one
element of infinite order in the basis.
\end{proof}

 Theorem
\ref{th:tipo} describes how to find the additive structure, in general, for any
 indecomposable $R$-module after
computing the $p$-height of the elements that connect the building
blocks in $M$.
Let $d_{ji}$ denote $d(j,i)$. Also let
 $\texttt{d}_i=p_2^{d_{2i}-1}s_{2i}=-p_1^{d_{1i+1}-1}s_{1i+1}$
for $1\leq i\leq m-1$. If $d_{11}\neq \infty$, let $\texttt{d}
_0=p_1^{d_{11}-1}s_{11}$,
 and  if $d_{2m}\neq \infty$, let $\texttt{d}_m=p_2^{d_{2m}-1}s_{2m}$.
The $p$-height of the element $\texttt{d}_i=p^{k-1}d'$ is  $h_p(\texttt{d}_i)=k-1$  in the Abelian group $t(M)$.
Let $\ell _\alpha $  be the number of elements $\texttt{d}_i$ such that
$h_p(\texttt{d}_i)=\alpha -1$, see \cite{F}.

  Let  $S=\bigoplus _i S^{(i)}$ be  a separated representation of an
 indecomposable $R$-module $M$. If $M$ is \emph{block cyclic},  then we
  consider the  separated module
 $S'=\oplus_iS'^{(i)}$ such that  $d'_{ji}=d_{ji}-1$ for all $(j,i)$. If $M$ is
 \emph{deleted cyclic}, then we consider the  separated module
 $S'=\oplus_iS'^{(i)}$ such that  $d'_{ji}=d_{ji}-1$ for  $(j,i)\ne (1,1)$
 and  $(j,i)\ne (2,m)$.

 \begin{theorem}\label{th:tipo}
 Let $M$ be an indecomposable $R$-module, and let $S=\bigoplus
 _{i=1}^m S^{(i)}$ be the separated representation of $M$. Then
 \[ \underline t(M)= \sum_{i=1}^{m}\underline t(S'^{(i)}) + (0,\ell
 _1-\ell _2,\ldots , \ell _{n-1}-\ell _n,\ell_n),\]
 where $n=\ex (t(M))$.
 \end{theorem}

 \begin{proof}
First, suppose that $M=\la a_1,\ldots ,a_m\ra $ is a deleted cycle
indecomposable $R$-module. If $m=1$, the theorem is obviously
true. Suppose the result is proved for $m-1$. Then consider $M/\la
\texttt{d}_1\ra= S'^{(1)}\oplus M'$, where $M'=\la a_2,\ldots
,a_n\ra /\la \texttt{d}_1\ra $. Since $M'$ is generated by $m-1$
elements, we can apply induction.  By \cite[Corollary 3.3]{AB},
if $h_p(\texttt{d}_1)=\beta$ then
\[ \begin{array}{ll}\underline t(M)&=\underline t(M/\la
\texttt{d}_1\ra)+v(\beta)
\\
&= \underline t (S'^{(1)})+\sum_{k=2}^n \underline t(S'^{(k)})+
(0,\ell ' _1-\ell ' _2,\ldots , \ell '_{n-1}-\ell ' _n,\ell '_n)+
v(\beta),
\end{array}\]
where $\ell '_\alpha $ is  the number of elements
$\texttt{d}'_k$ such that $h_p(\texttt{d}'_k)=\alpha -1$ in $M'$.
The vector $v(\beta)=(v_0,v_1, \ldots ,v_n)\in \Z^{n+1}$ satisfies
$v_{\beta -1}=-1$, $v_\beta =1$ and $v_i=0$ otherwise. It is
clear that $h_p(\texttt{d}'_k)=h_p(\texttt{d}_k)$ for $k\geq 2$.
So our claim holds.

Now suppose $M$ is a block cyclic indecomposable module. Observe that
the result holds for the module $M'=M/\la \texttt{d}_0\ra$, since $M'$ is a deleted
cyclic module. Hence our claim holds by \cite[Corollary 3.3]{AB}.
 \end{proof}


\begin{thebibliography}{99}

\bibitem {AL}  D. Arnold, R. Laubenbacher, {\it Finitely generated modules over pullback
rings.}, J. Algebra \textbf{184} (1996), pp. 304--332.


\bibitem{A} Michael Artin, {\it Algebra }, Prentice Hall, 1991.

\bibitem{H} Marshall Hall, {\it The theory of Groups}, The Macmillan Company, New
York, 1959.

\bibitem{AB} M.  A. Avi\~n\'o  and R. Bautista, {\it The Additive Structure of
Indecomposable $\Z_{p^n}C_p$-Modules}, Communications in Algebra.
{\bf 24} (1996), no. 8.

\bibitem{Bu} B. Buchberger, {\it Gr\"obner bases-an algorithmic method in polynomial
ideal theory,} Chapter 6 in N.K. Bose (ed.): Multidimensional
Systems Theory, D. Reidel Publ., (1985).

\bibitem{CLO} D. Cox J. Little, D. O'Shea, {\it Ideals, Varieties, and Algorithms},
Springer, Second Edition, 1996.

\bibitem{CR} C. W.  Curtis and I. Reiner, {\it Representation Theory
of Finite groups and Associative Algebras}, John Wiley and Sons,
1962.

\bibitem{KM} M. I. Kargapolov, Ju. I. Merzjgakov, {\it Fundamentals
of the Theory of Groups}, Springer-Verlag, New York INc, 1979.


\bibitem{F}  L. Fuchs, {\it Infinite Abelian Groups},  Vol. \textbf{I/II},  Academic Press, 1970/73.

\bibitem{GP} I. M. Gelfand , V. A. Ponomariev, {\it Indecomposable
representations of a Lorentz group} Uspekhi Matem. Nauk, 140, 3-60
(1968).

\bibitem {G} Daniel Gorenstein, {\it Finite Groups}, Harper's Series
in Modern Mathematics, 1968.

\bibitem{L1} L. S. Levy, {\it Modules over Pullbacks and Subdirect Sums}, J. Algebra {\bf 71} (1981), 50-61.


\bibitem{L2} L. S. Levy, {\it Mixed Modules over $\Z G$, $G$ Cyclic of
Prime Order, and over Related Dedekind Pullbacks}, J. Algebra {\bf
71} (1981), 62-114.


\bibitem{NR} L. A. Nazarova,  A. V.  Roiter, {\it Finitely generated
modules over a Dyad of two local Dedekind rings, and finite groups
with an abelian normal divisor of index $p$, } Izv. Akad. Nauk.
SSSr. Ser. Mat. Tom 33, (1969), No 1, 65-86.

\bibitem{NRSB}L. A. Nazarova, A. V. Roiter,V. V.  Sergeitchuck, V. M.  Bondarenko, {\it
Application of Modules over a Dyad for the classification of finite
$p$-groups  possessing an abelian subgroup of index $p$ and of pairs
of mutually annihilating operators.} V. A. Steklova Akademii Nauk
SSSR, Vol. 28, pp 69-92, 1972.


\bibitem{S} B. Sturmfels, {\it Gr\"obner Bases and Convex Polytopes}. University Lectures
Series. Vol {\bf 8}, American Mathematical Society, Providence,
1996.

\bibitem{Sz} G. Szekeres, {\it Determination of Certain Family of Finite Metabelian Groups},
 Trans. Amer. Math. Soc., {\bf 66}, 1949, 1-43.
 \end{thebibliography}
 \end{document}